\documentclass[11pt,amsart]{article}
\usepackage[all]{xy}
\usepackage[centertags]{amsmath}
\usepackage{latexsym,amsthm,amssymb,amsmath,verbatim}
\usepackage{epsf}
\newtheorem{thm}{Theorem}[section]
\newtheorem{prop}[thm]{Proposition}
\newtheorem{lem}[thm]{Lemma}

\newtheorem{cor}[thm]{Corollary}
\newtheorem{defn}[thm]{Definition}
\newtheorem*{defn*}{Definition}
\newtheorem{rem}[thm]{Remark}

\newcommand{\des}{{\rm des}}
\newcommand{\maj}{{\rm maj}}

\newcommand{\Colrm}{{\rm Col}}
\newcommand{\col}{{\rm col}}

\newcommand{\Des}{{\rm Des}}

\newcommand{\la}{\lambda}
\newcommand{\s}{\sigma}
\newcommand{\g}{\gamma}

\newcommand{\x}{{\bf x}}
\newcommand{\z}{{\bf z}}
\newcommand{\y}{{\bf y}}

\newcommand{\G}{G(r,p,n)}
\newcommand{\bde}{\begin{defn}}
\newcommand{\ede}{\end{defn}}
\newcommand{\fmaj}{\operatorname{fmaj}}
\newcommand{\fdes}{\operatorname{fdes}}

\newcommand{\inv}{\operatorname{inv}}

\newtheorem{exa}[thm]{Example}
\newcommand{\h}{\hat{g}}
\newcommand{\lab}{{\vec{\lambda}}}
\newcommand{\alp}{{\vec{\alpha}}}
\newcommand{\muv}{{\vec{\mu}}}

\newcommand{\T}{{\mathcal{T}}}

\newcommand{\GA}{\Gamma}
\newcommand{\grn}{G(r,n)}

\newcommand{\PP}{\mathbb{P}}
\newcommand{\NN}{\mathbb{N}}
\newcommand{\ZZ}{\mathbb{Z}}
\newcommand{\CC}{\mathbb{C}}

\newcommand{\sumlim}{\sum\limits}
\newcommand{\prodlim}{\prod\limits}
\newcommand{\Sno}{{S_{012345}}}
\newcommand{\Sone}{{S_{234501}}}
\newcommand{\Stwo}{{S_{450123}}}

\begin{document}

\title{Colored-Descent Representations of Complex Reflection Groups $G(r,p,n)$}

\author{Eli Bagno and Riccardo Biagioli}

\date{\today}

\maketitle
\begin{abstract}
We study the complex reflection groups $G(r,p,n)$. By considering
these groups as subgroups of the wreath products $\mathbb{Z}_r \wr
S_n$, and by using Clifford theory, we define combinatorial
parameters and descent representations of $G(r,p,n)$, previously
known for classical Weyl groups. One of these parameters is the
flag major index, which also has an important role in the
decomposition of these representations into irreducibles. A
Carlitz type identity relating the combinatorial parameters with
the degrees of the group, is presented.
\end{abstract}

\section{Introduction}

Let $V$ be a complex vector space of dimension $n$.  A {\em
pseudo-reflection} on $V$ is a linear transformation on $V$ of
finite order which fixes a hyperplane in $V$ pointwise. A {\em
complex reflection group} on $V$ is a finite subgroup $W$ of ${\rm
GL}(V)$ generated by pseudo-reflections. Such groups are
characterized by the structure of their invariant ring. More
precisely, let $\CC[V]$ be the symmetric algebra of $V$ and let us
denote by $\CC[V]^W$ the algebra of invariants of $W$. Then
Shephard-Todd \cite{ST} and Chevalley \cite{Ch} proved that $W$ is
generated by pseudo-reflections if and only if $\CC[V]^W$ is a
polynomial ring.

Irreducible finite complex reflection groups have been classified
by Shephard-Todd \cite{ST}. In particular, there is a single
infinite family of groups and exactly 34 other ``exceptional"
complex reflection groups. The infinite family $\G$, where $r,p,n$
are positive integers numbers with $p | r$, consists of the groups
of $n\times n$ matrices such that:

1) the entries are either 0 or $r^{\rm th}$ roots of unity;

2) there is exactly one nonzero entry in each row and each column;

3) the $(r/p)^{{\rm th}}$ power of the product of the nonzero entries is 1.

\noindent The classical Weyl groups appear as special cases:
$G(1,1,n)=S_n$ the symmetric group, $G(2,1,n)=B_n$ the hyperoctahedral group, and $G(2,2,n)=D_n$ the group of even-signed permutations.

Through research on complex reflection groups and their braid
groups and Hecke algebras, the fact that complex reflection groups
behave like Weyl groups has become more and more clear. In
particular, it has been recently discovered that they (and not only
Weyl groups) play a key role in the structure as well as in the
representation theory of finite reductive groups. For more
information on these results the reader is advised to consult the
survey article of Brou\'e \cite{B}, and the handbook of Geck and
Malle \cite{GM}.

One of the aims of this paper is to show that complex reflection
groups continue to behave like (classical) Weyl groups also from
the point of view of the combinatorial representation theory. It
is well known that, in a way similar to Coxeter groups, they have
presentations in terms of generators and relations, that can be
visualized by Dynkin type diagrams (see e.g., \cite{BMR}).
Moreover, their elements can be represented as colored
permutations. In fact, the complex reflection groups $\G$ can be
naturally identified as normal subgroups of index $p$ of the
wreath product $\grn:=\ZZ_r \wr S_n$, where $\mathbb{Z}_r$ is the
cyclic group of order $r$. This makes it possible to handle them
by purely combinatorial methods. We follow this approach.

The organization of the paper is as follows. In Section
\ref{flag}, we introduce several combinatorial parameters on
complex reflection groups. Among those we define the concept of
``major index" and ``descent number" for $\G$. Then our
investigation continues by showing the interplay between these new
combinatorial statistics and the representation theory of $\G$.
More precisely, in Section \ref{col-des} we define a new set of
$G(r,p,n)$-modules, which we call {\em colored-descent
representations}. They generalizes to all groups $G(r,p,n)$, the
descent representations introduced for $S_n$ and $B_n$, by Adin,
Brenti, and Roichman \cite{ABR}. These modules are isomorphic to
the Solomon's descent representations in the case of $S_n$, while
they decompose Solomon's representations in the case of $B_n$.
Following Adin-Brenti-Roichman's approach, we use the coinvariant
algebra as representation space. In order to do that, in Section
\ref{col-des-basis}, we use refined statistics to define an
explicit monomial basis for the coinvariant space. This basis has
special properties that allow us to define our new set of modules.
In Section \ref{S-decomposition}, the decomposition into
irreducibles of the colored-descent representations is provided
(Theorem \ref{main}). We use a generalization of a formula of
Stanley (Theorem \ref{4.5}) on a specialization of a product of
Schur functions, proved in Section \ref{S-stanley}. It turns out
that the multiplicity of any irreducible representations is
counted by the cardinality of a particular class of standard Young
tableaux, called $n$-{\em orbital}. As a corollary of that, we
obtain a refinement of a theorem first attributed to Stembridge
(Corollary \ref{stembrH}). Finally, in Section \ref{identities},
by using the properties of the colored descent basis, a Carlitz
type identity (Theorem \ref{Ca-H}) relating the major index and
descent number with the degrees of $G(r,p,n)$ will be derived.

It is worth to say, that in a very recent preprint concerning wreath products, Baumann and Hohlweg \cite{BH}, following Solomon's descent algebra approach
 \cite{So},
 define some special characters of $G(r,n)$ as images
 through a ``generalized Solomon homomorphism" of the elements of a  basis of the Mantaci-Reutenaur algebra \cite{MR}. Nevertheless, they don't provide any module having them as characters. It was a nice surprise to see that the modules
which carry them as characters are our colored-descent representations
in the case of $G(r,1,n)$.

\section{Complex Reflection Groups $\G$}\label{generators}

For our exposition it will be much more convenient to consider
wreath products not as groups of complex matrices, but as groups
of colored permutations.

Let $\mathbb{P}:=\{1,2,\ldots\}$, $\NN := \PP \cup \{ 0 \}$,
 and $\CC$ be the field of complex numbers. For any $n\in \PP$,  we let
$[n]:=\{1,2,\ldots,n\},$ and for any $a, b \in \NN$ we let
$[a,b]:=\{a,a+1,\ldots,b\}.$ Let $S_n$ be the symmetric group on
$[n]$. A permutation $\s \in S_n$ will be denoted by $\s=\s(1)\cdots
\s(n).$

Let $r,n \in \mathbb{P}$. The {\em wreath product} $\grn$ of
$\mathbb{Z}_r$ by $S_n$ is defined by
\begin{equation}\label{def-grn}
G(r,n):=\{((c_1,\ldots,c_n),\sigma) \mid c_i \in [0,r-1],\sigma \in
S_n\}.
\end{equation}
Any $c_i$ can be considered as the color
of the corresponding entry $\s(i)$. This explains the fact that
this group is also called the {\em group of r-colored
permutations}. Sometimes we will represent its elements in {\em window
notation} as
\[g=g(1)\cdots g(n)=\s(1)^{c_1}\cdots\s(n)^{c_n}.\]
When it is not clear from the context, we will denote $c_i$ by
$c_i(g)$. Moreover, if $c_i=0$, it will be omitted in the window notation of
$g$. We denote by
\[\Colrm(g):=(c_1,\ldots,c_n) \;\;\; {\rm and} \; \;\; \col(g): =\sum_{i=1}^n c_i,\]
the {\em color vector} and the {\em color weight} of any $g:=((c_1,\ldots,c_n),\s) \in G(r,n)$.\\
\noindent For example, for $g=4^1 3 2^4 1^2 \in G(5,4)$ we have
$\Colrm(g)=(1,0,4,2)$ and $\col(g)=7$.
\bigskip

Now let $p \in \mathbb{P}$ be such that $p|r$. The {\em complex
reflection group} $G(r,p,n)$ is the subgroup of $G(r,n)$ defined
by
\begin{equation}\label{def-grpn}
G(r,p,n):=\{g \in G(r,n) \mid \col(g)\equiv 0  \; {\rm mod} \; p\}.
\end{equation}
Note that $\G$ is the kernel of the map $G(r,n) \longrightarrow \ZZ_p$, sending
$g$ in its color weight $\col(g)$, and so it is a normal subgroup of $G(r,n)$ of index $p$.
It is clear from the definition that the wreath product is $G(r,1,n)$. Moreover, $G(1,1,n)$ is the symmetric group, $G(2,2,n)$ is the Weyl groups of type $D$, and $G(r,r,2)$
is the dihedral group of order $2r$.

\section{Flag Major Index on $G(r,p,n)$} \label{flag}

The major index is a very well studied statistic over the
symmetric group, (see e.g., \cite{Ca},\cite{Fo},\cite{GG},\cite{Ge},\cite{MM}) .
Lately, many generalizations of this concept have been given for
the other classical Weyl groups and for wreath products (see e.g., \cite{AR},\cite{ABR1},\cite{BC},\cite{Rein},\cite{Steing}). Among those statistics, it became
soon clear that, the flag major index,
introduced by Adin and Roichman on $B_n$ \cite{AR}, is the one more similar to the major index \cite{GG}, regarding the aspect of its implications with the representation theory of the group \cite{AR,ABR}. Its analogue for Weyl groups of type $D$ is the
$D$-(flag) major index, introduced in \cite{BC}. In this section we
define a version of the major index for $\G$. We introduce a
particular subset of $G(r,n)$ which is in bijective correspondence with $\G$, and
define there our statistics. In the following sections, it will
become clear that this subset is the right object to work with.

\bigskip
In order to lighten the notation, we let $G:=\grn$ and $H:=\G$. For any $r,p,n \in \mathbb{P}$, with $p|r$ let $d:=r/p$. We
define the following subset of $G(r,n)$,
\begin{equation}\label{def-Gamma}
\GA(r,p,n)=\{\g=((c_1,\ldots,c_n),\s) \in G(r,n) \mid c_n < d\}.
\end{equation}
Note that $\GA:=\Gamma(r,p,n)$ is not a subgroup of $G$. Clearly,
$|\Gamma|=n! r^{n-1}d$ and so it is in bijection with $H$.
Moreover, one can easily check that the mapping $\varphi$
\begin{equation}\label{bijection}
((c_1,\ldots,c_n),\s) \stackrel{\varphi}{\mapsto} ((c_1,\ldots ,\lfloor \frac{c_n}{p} \rfloor),\s)
\end{equation}
is a bijection between $H$ and $\Gamma$. As usual, for any $a \in
\mathbb{Q}$, $\lfloor a \rfloor$ denotes the greatest integer
$\leq a$.

In order to make our definitions more natural and clear, from now
on, we will work with $\GA$ instead of $H$. Clearly, via the above
bijection $\varphi$ every function on $\GA$ can be considered as a function
on $H$ and viceversa.

\bigskip

We fix the following order $\prec$ on colored integer
numbers
\begin{equation}\label{order}
1^{r-1}\prec 2^{r-1}\prec \ldots \prec n^{r-1} \prec \ldots \prec 1^{1} \prec 2^{1}\prec \ldots \prec n^{1} \prec 1 \prec 2 \prec \ldots \prec n.
\end{equation}
The {\em descent set} of a colored integer sequence $\g \in \GA$
is defined by
\[\Des(\g):=\{i \in [n-1] : \g_i \succ \g_{i+1} \}.\]
The {\em major index} of $\g$ is the sum of all the descents of $\g$, i.e.,
\[\maj(\g):=\sum_{i \in \Des(\g)} i.\]

Following \cite{AR}, we define the {\em flag major index} of $\g \in \Gamma$ by
\[\fmaj(\g):=r\cdot \maj(\g) + \col(\g).\]
Via the identification given by the bijection $\varphi$ we
consider this as the flag major index on complex reflection
groups. Often for our purposes, it will be better to obtain this
statistic in a more elaborate way. For any
$\g=((c_1,\ldots,c_n),\s) \in \GA$ we let
\begin{eqnarray}\label{def-icomp}
d_i(\g)&:=&|\{j \in \Des(\g): j \succeq i\}| \;\; {\rm and}\\
f_i(\g)&:=&r \cdot d_i(\g) + c_i(\g).
\end{eqnarray}
For any $\g \in \Gamma$, we define the {\em flag descent number} of $\g$ by
\begin{equation}\label{def-fdes}
\fdes(\g):=r\cdot d_1(\g) + c_1.
\end{equation}
It is clear that for every $\g \in \GA$,
$\fmaj(\g):=\sumlim_{i=1}^n f_i(\g)$, and that $d_1(\g)$ is the cardinality of the set $\Des(\g)$. Let us conclude this section with an example.
\begin{exa}\label{ex-1}
{\rm Let $\g={\bf 6}2^{5}4^{4}{\bf 3}^{1}1^{6}5^{3} \in \Gamma(8,2,6)$. The set of descents is $\Des(\g)=\{1,4\}$. Hence $(d_1(\g),\ldots,d_n(\g))=(2,1,1,1,0,0)$,
$(f_1(\g),\ldots,f_n(\g))=(16,13,12,9,6,3)$, and so $\fmaj(\g)=59$ and $\fdes(\g)=8$.}
\end{exa}

\section{Colored Descent Basis}\label{col-des-basis}

Let $W \leq {\rm GL(V)}$ be a complex reflection groups. If we set $\x=x_1,\ldots,x_n$ as a basis for $V$, then $\CC[V]$ can be identified with the ring of polynomials
$\CC[\x]$. The ring of invariants $\CC[\x]^W$ is then generated by
1 and by a set of $n$ algebraically independent homogeneous
polynomials $\{\vartheta_1,\ldots , \vartheta_n\}$ which are
called {\em basic invariants}. Although these polynomials are not
uniquely determined, their degrees $d_1,\ldots,d_n$ are basic
numerical invariants of the group, and are called the {\em degrees
of} $W$. Let us denote by $\mathcal{I}_W$ the ideal generated by
the invariants of strictly positive degree. The {\em
coinvariant space} of $W$ is defined by
\[\CC[\x]_W:=\CC[\x]/\mathcal{I}_W.\]
Since $\mathcal{I}_W$ is $W$-invariant, the group $W$ acts naturally on $\CC[\x]_W$. It is well known that $\CC[\x]_W$ is isomorphic to the {\em
left regular representation} of $W$. It follows that its dimension
 as a $\CC$-module is equal to the order of the
group $W$.

As a first application of the flag major index, we find a basis of
the coinvariant space of $H$. The wreath product $G$ acts on the ring of
polynomials $\CC[\x]$ as follows
\begin{equation}\label{def-action}
\s(1)^{c_1}\cdots\s(n)^{c_n} \cdot P(x_1,\ldots,x_n)=P(\zeta^{c_{\s(1)}}x_{\s(1)},\ldots,\zeta^{c_{\s(n)}}x_{\s(n)}),
\end{equation}
where $\zeta$ denotes a primitive $r^{\rm th}$ root of unity.
A set of basic invariants under this actions is given by the
{\em elementary symmetric functions} $e_j(x_1^r,\ldots,x_n^r)$,
$1\leq j \leq n$. Now, consider the restriction of the previous
action on $\CC[\x]$ to $H$. A set of fundamental invariants is given by
\[\vartheta_j(x_1,\ldots,x_n):=\left\{ \begin{array}{ll} e_j(x_1^r,\ldots,x_n^r) & {\rm for} \;  j=1,\ldots,n-1
\\
x_1^{d}\cdots x_n^{d} & {\rm for} \; j=n.\end{array} \right.\] It
follows that the degrees of $H$ are $r,2r,\ldots,(n-1)r,nd$. Let $\mathcal{I}_H:=(\vartheta_1,\ldots,\vartheta_n)$ be the ideal generated by the constant term invariant polynomials of $H$. The coinvariant space, $\CC[\x]_H:=\CC[\x]/ \mathcal{I}_H$, has dimension equal to $|H|$, that is $n! r^n/p$. In what follows we will associate to any element $h \in H$ an ad-hoc monomial in $\CC[\x]$. Those
monomials will form a linear basis of  $\CC[\x]_H$.
\bigskip

Let $\g=((c_1,\ldots,c_n),\s)\in \GA$. We define
\begin{equation}\label{defmonomial}
\x_{\g}:=\prodlim_{i=1}^n x_{\s(i)}^{f_i(\g)}.
\end{equation}
By definition $f_n(\g)<d$, hence $\x_{\g}$ is nonzero in
$\CC[\x]_H$. Clearly ${\rm deg} (\x_\g)=\fmaj(\g)$.\\
\noindent For example, let $\g=62^{5}4^{4}3^{1}1^{6}5^{3} \in
\Gamma(8,2,6)$ as in Example \ref{ex-1}. The associated monomial is $\x_{\g}=x_1^{6} x_2^{13} x_3^9 x_4^{12}
x_5^{3} x_6^{16}.$

\bigskip

We restrict our attention to the quotient $
S:=\CC[\x]/(\vartheta_{n}) $. Hence we consider nonzero monomials
$ M=\prodlim ^{n}_{i=1}x^{a_{i}}_{i} $  such that $a_{i} < d$ for
at least one $i \in [n]$. We associate to $M$ the unique element
$\gamma (M)=((c_1,\ldots,c_n),\s) \in \GA$ such that for all $i
\in [n]$:
\begin{itemize}
\item [$i)$] $ a_{\s (i)}\geq a_{\s(i+1)} $;
\item [$ii)$] $ a_{\s (i)}=a_{\s (i+1)}\Longrightarrow \s (i) < \s (i+1) $,
\item [$iii)$] $a_{\s (i)} \equiv  c_i \: ({\rm mod} \: r)$.
\end{itemize}
We denote by $\la(M):=(a_{\s (1)},\ldots ,a_{\s (n)})$ the
\emph{exponent partition} of $M$, and we call $ \g(M) \in \GA $
the \emph{colored index permutation}.

Now, let $M=\prodlim ^{n}_{i=1}x^{a_{i}}_{i}$ be a nonzero
monomial in $S$, and let $\g:=\gamma (M)$ be its colored
index permutation. Consider now the monomial $\x_{\g}$ associated to
$\g$. It is not hard to see that the sequence $ (a_{\s
(i)}-f_{i}(\gamma )),\, i=1,\ldots ,n-1 $, of exponents of
$M/\x_{\g}$, consists of nonnegative integers of the form $c\cdot
r$ with $c>0$, and is weakly decreasing.
This allows us to associate to $M$ the {\em complementary partition} $\mu(M)$,
defined by
\begin{equation}\label{def-mu}
\mu^{\prime}(M):=\left(\frac{a_{\s (i)}-f_{i}(\gamma )}{r}\right)^{n-1}_{i=1},
\end{equation}
where, as usual, $\mu^{\prime}$ denotes the conjugate partition of
$\mu$.


\begin{exa}{\rm
Let $r=8$, $p=2$, and $n=6$ and consider the monomial
$M=x^{6}_{1}x_{2}^{21}x^{17}_{3}x_4^{20}x_{5}^{3}x^{32}_{6} \in
\CC[x_1,\ldots,x_6]/(\vartheta_6)$. The exponent partition $\la(M)=(32,21,20,17,6,3)$ is obtained by
reordering the power of $x_i$'s following the colored index permutation $\gamma
(M)=62^{5}4^{4}3^{1}1^{6}5^{3} \in \Gamma(8,2,6)$. We have already
computed the monomial
$\x_{\g(M)}=x_1^6x_2^{13}x_3^9x_4^{12}x_5^3x_6^{16}$. It follows
that $\mu (M)=(4,1).$}
\end{exa}

We now define a partial order $\sqsubseteq$ on the monomials of the same total
degree in $ S $. Let $ M $ and $ M^{\prime } $ be nonzero monomials in $ S $
with the same total degree and such that the exponents of $ x_{i}
$ in $ M $ and $ M^{\prime } $ have the same parity (mod $r$)
for every $ i\in [n] $. Then we write $ M^{\prime } \sqsubset M $ if one
of the following holds:
\begin{itemize}
\item [$1)$] $ \lambda (M^{\prime })\lhd \lambda (M) $, or
\item [$2)$] $ \lambda (M^{\prime })=\lambda (M)\, \textrm{ and }\inv(\gamma (M^{\prime }))>\inv(\gamma (M)) $.
\end{itemize}
Here, $\inv(\g):=|\{(i,j) \; | \; i<j \;\; {\rm and} \;\; \g(i)
\succ \g(j)\}|,$ and  $\lhd$ denotes the \emph{dominance order}
defined on the set partitions of a fixed nonnegative integer $ n
$ by:  $ \mu \unlhd \lambda  $ if for all $ i\geq 1 $
\[ \mu _{1}+\mu_{2}+\cdots +\mu _{i}\leq \lambda _{1}+\lambda _{2}+\cdots+\lambda _{i}.\]

By arguments similar to those given in the proof of Lemma 3.3 in \cite{BC2},
it can be showed that every monomial $M \in S$ admits the following
 expansion in terms of the basis monomials $\x_{\g}$'s and the basic generators
 $\{\vartheta_i\}$ of the ideal $\mathcal{I}_H$:
\begin{equation}\label{straightening}
M=\vartheta_{\mu (M)}\cdot \x_{\gamma (M)}+\sum _{M^{\prime
} \sqsubset M}n_{M^{\prime },M}\vartheta_{\mu (M^{\prime })}\cdot \x_{\gamma
(M^{\prime })},
\end{equation}
where $ n_{M,M^{\prime }} $ are integers. As usual,
$\vartheta_{\mu}:=\vartheta_{\mu_1}\vartheta_{\mu_2}\cdots \vartheta_{\mu_{\ell}},$ with
$\ell:=\ell(\mu)$ the {\em length} of the partition $\mu$.

\begin{exa}\label{ex-M}{\rm
Consider the group $ \Gamma(6,2,3)$ and the monomial $
M=x^{11}_{1} x_{2}^8 x_{3}. $ We have $ \gamma
(M)=1^{5}2^{2}3^{1} \in \Gamma(6,2,3)$, $ \x_{\gamma
(M)}=x^{5}_{1} x_{2}^2 x_{3} $, and $ \mu (M)=(2) $. Then, if
we set $ M_{1}=x^{11}_{1}x^{2}_{2}x^{7}_{3} $, and $ M_{2}=x^{5}_{1}x^{8}_{2}x^{7}_{3} $ we have that
$M=\x_{\gamma (M)}\vartheta_{2}-M_{1}-M_{2}$. It is easy to see that  $ M_{1}, M_{2} \sqsubset M $ and
that $M_1=\x_{\gamma (M_1)}$, and $M_2=\x_{\gamma (M_2)},$ so
\[M=\x_{\gamma (M)}\vartheta_{2}-\x_{\g(M_1)}-\x_{\g(M_2)}.\]}
\end{exa}

Since the monomials $\{\x_{\g} + \mathcal{I}_H \mid \g \in \Gamma
\}$ are generators for $\CC[\x]_H$, and ${\rm dim} \; \CC[\x]_H=|H|$, they
form a basis, called the {\em colored-descent basis}. We summarize
this by:

\begin{thm}\label{basis-coivariants}
The set
\[\{\x_{\g} + \mathcal{I}_H : \; \g \in \GA \}\]
is a basis for $\CC[\x]_H$.
\end{thm}

Note that, when $H$ specializes to one of the classical Weyl groups,
our basis coincides with the descent basis defined by Garsia-Stanton  for $S_n$ \cite{GS},
by Adin-Brenti-Roichman  for $B_n$ \cite{ABR},
and by Biagioli-Caselli for $D_n$ \cite{BC2}.
Recently, another basis for $\CC[\x]_H$ has been given by Allen \cite{Al}.
Although both our and Allen's basis
coincide with the Garsia-Stanton basis in the case of $S_n$, in general they are different as
can be checked already in the small case of $G(2,2,2)$.
It would be interesting to see if  Allen's basis leads to
an analogous definition of descent representations (see Section \ref{col-des}).

\begin{exa}{\rm
The elements of $\Gamma(6,3,2)$, $d=2$, are
\[\begin{array}{llllll}
12 & 1^{1}2 &1^{2}2 &1^{3}2 &1^{4}2 &1^{5}2 \\
12^{1} & 1^{1}2^{1} &1^{2}2^{1} &1^{3}2^{1} &1^{4}2^{1} &1^{5}2^{1} \\
21 & 2^{1}1 & 2^{2}1 & 2^{3}2 & 2^{4}1 & 2^{5}1 \\
21^{1} & 2^{1}1^{1} & 2^{2}1^{1} & 2^{3}2^{1} & 2^{4}1^{1} & 2^{5}1^{1}.
\end{array}\]
The corresponding monomials
\[\begin{array}{llllll}
1 & x_1 & x_1^2 & x_1^3 & x_1^4 & x_1^5 \\
x_1^6 x_2 & x_1 x_2 & x_1^2 x_2 & x_1^3 x_2 & x_1^4 x_2 & x_1^5 x_2 \\
x_2^6 & x_2 & x_2^2 & x_2^3 & x_2^4 & x_2^5 \\
x_2^6 x_1 & x_2^7 x_1 & x_2^2 x_1 & x_2^3 x_1 & x_2^4 x_1 & x_2^5 x_1\\
\end{array}\]
form a basis for $\CC[x_1,x_2]/(x_1^6+x_2^6, x_1^2x_2^2)$.}
\end{exa}

%
%

\section{The Representation Theory of $\G$}\label{repr-theory}

In this section we present the representation theory of the group
$H:=\G$. We follow the exposition of \cite{Ste}, (see also \cite{MY}). Since the
irreducible representations of $H$ are related to the irreducible
representations of $G$ via Clifford Theory, we start this
section by presenting the representation theory of $G$.

Let  $g=\s(1)^{c_1}\cdots\s(n)^{c_n} \in G$. First divide $\s \in
S_n$ into cycles, and then provide the entries with their original
colors $c_i$, thus obtaining {\it colored cycles}. The color of a
cycle is simply the sum of all the colors of its entries.
%
For every $i \in [0,r-1]$, let $\alpha^i$ be the partition formed
by the lengths of the cycles of $g$ having color $i$. We may thus
associate $g$ with the $r$-partition
$\alp=(\alpha^0,\ldots,\alpha^{r-1})$. Note that
$\sumlim_{i=0}^{r-1}{|\alpha^i|}=n$. We refer to $\alp$ as the
{\em type} of $g$. For example, the decompositions in colored cycles of $\g=62^{5}4^{4}3^{1}1^{6}5^{3} \in G(8,6)$ is $(1^665^3)(3^14^4)(2^5)$. We have $\alpha^1=(3)$, $\alpha^5=(2,1)$, and $\alpha^i=0$ for all other $i$'s.

One can prove that two elements of $G$ are
conjugate if and only if they have the same type. It is well known
that irreducible representations of $G$ are also indexed by
$r$-tuple of partitions $\lab:=(\la^0,\ldots,\la^{r-1})$ with
$\sumlim_{i=0}^{r-1}{|\lambda^i|}=n$. We denote this set by
$\mathcal{P}_{r,n}$.

As mentioned above, the passage to the representation theory of
$H$, is by Clifford theory. The group $G/H$ can be identified with
the cyclic group $C$ of order $p$ consisting of the characters
$\delta$ of $G$ satisfying $H \subset {\rm Ker}(\delta)$. More
precisely, define the linear character $\delta_0$ of $G$ by
$\delta_0((c_1,\ldots,c_n),\s):=\zeta^{c_1+\ldots+c_n}$, so that
$C=<\delta_0^d>\simeq \mathbb{Z}_p$.

The group $C$ acts on the set of irreducible representations of
$G$ by $V(\lab) \mapsto \delta \otimes V(\lab)$, where $V(\lab)$
is the irreducible representation of $G$ indexed by $\lab$, and
$\delta \in C$. This action can be explicitly described as
follows. Let $\lab=(\lambda^0,\ldots,\lambda^{r-1}) \in
\mathcal{P}_{r,n}$. We define a {\em 1-shift} of $\lab$ by
\begin{equation}\label{shift-op}
(\lab)^{\circlearrowleft 1}:= (\la^{r-1},\la^0,\ldots,\la^{r-2}).
\end{equation}
By applying $i$-times the {\em shift operator} we get
$(\lab)^{\circlearrowleft i}$. Then one can show (see
\cite[Section 4]{MY}) that
\begin{equation}\label{action}
\delta_0 \otimes V(\lab) \simeq V((\lab)^{\circlearrowleft 1}),
\end{equation}
for every $\vec{\lambda} \in \mathcal{P}_{r,n}$. Now let us denote by $[\lab]$ a $C$-orbit of the representation
$V(\lab)$. From (\ref{action}) we obtain that
$[\lab]=\{V(\vec{\mu}) \; : \; \vec{\mu} \sim \lab\}$, where the
equivalence relation is defined by
\begin{equation}\label{orbit}
\lab \sim \vec{\mu} \;\; {\rm  if \; and \; only} \;\;  \vec{\mu}
= (\lab)^{\circlearrowleft i \cdot d}\;\; {\rm for}  \;\; {\rm
some}\;\; i \in [0,p-1].
\end{equation}
Let us denote $b(\lab):=|[\lab]|$, and set
$u(\lab):=\frac{p}{b(\lab)}$. Consider the stabilizer of $\lab$,
$C_{\lab}$, that is:
$$C_{\lab}:=\{\delta \in C \mid V(\lab)=\delta \otimes V(\lab)\}.$$
Clearly, $C_{\lab}$ is a subgroup of $C$ generated by
$\delta_0^{b(\lab) \cdot d}$ and so $|C_{\lab}|=u(\lab)$.

It can be proven that the restriction of the irreducible
representation $V(\lab)$ of $G$ to $H$, decomposes into
$u(\lab)=|C_{\lab}|$ non-isomorphic irreducible $H$ modules. On
the other hand, any other $G$-module in the same orbit $[\lab]$
will give us the same result. Actually, one can prove even more (see e.g., \cite{Ste}).

\begin{thm}\label{repsofgrpn}
There is a one to one correspondence between the irreducible
representations of $H$ and the ordered pairs $([\lab],\delta)$
where $[\lab]$ is the orbit of the irreducible representation
$V(\lab)$ of $G$ and $\delta \in C_{\lab}$. Moreover, if
$\chi^{\lab}\!\!\!\downarrow^G_H$ denotes the restriction of the character of
$V(\lab)$ to $H$, then
\begin{eqnarray*}
&i)& \chi^{\lab}\!\!\!\downarrow^G_H=\chi^{\vec{\mu}}\!\!\!\downarrow^G_H,\;\; for \; all \;\; \lab \sim \vec{\mu}, \; \; and  \\
&ii)& \chi^{\lab}\!\!\!\downarrow^G_H=\sumlim_{\delta \in
C_{\lab}}{\chi^{([\lab],\delta)}}.
\end{eqnarray*}

\end{thm}

Here is a simple but important example. The irreducible
representations of $B_n$ ($G(2,1,n)$ in our notation) are indexed
by bi-partitions of $n$. The Coxeter group $D_n$ ($G(2,2,n)$ in
our notation), is a subgroup of $B_n$ of index $2$. Thus the
stabilizer of the action of $B_n/D_n \cong \mathbb{Z}_2$ on a pair
of Young diagrams $\lab=(\lambda^1,\lambda^2)$ is either
$\mathbb{Z}_2$ if $\lambda^1=\lambda^2$, or $\{\rm id \}$ if
$\lambda^1 \neq \lambda^2$. In the first case, the irreducible
representation of $B_n$ corresponding to $\lab$, when restricted
to $D_n$, splits into two non-isomorphic irreducible
representations of $D_n$. In the second case
$\lab=(\lambda^1,\lambda^2)$ and $\lab^T=(\lambda^2,\lambda^1)$
correspond to two isomorphic irreducible representations of $D_n$.

\section{Frobenius Formula for $\G$}

Let $\y=y_1,y_2, \ldots, $ be an infinite set of variables and let
$\la$ be a partition of $n$. Let us denote by
$p_{\la}(\y):=p_{\la_1}(\y)p_{\la_2}(\y)\cdots p_{\la_{\ell}}(\y)$
the {\em power sum symmetric function}, where
\[p_k(\y):=y_1^k+y_2^k+\ldots .\]
A well know theorem of Frobenius establishes a link between the power sum and the Schur functions via the characters of the symmetric group. Let denote by $\chi^{\la}_{\alpha}$ the character of the irreducible
representation of $S_n$ corresponding to $\la$ calculated in the conjugacy class $\alpha$, and by $s_{\la}(\bf y)$ the {\em Schur function} corresponding to $\la$. Then  we have (see e.g., \cite[Corollary 7.17.4]{StaEC2})
\[ p_{\alpha }(\y)=\sum _{\lambda \vdash n}\chi
^{\lambda }_{\alpha }s_{\la }(\y).\] We need an analogous formula
involving the characters of the group $H$. In developing such an
analogue, we follow \cite[Appendix B]{MD} but we use a slightly
different notation.

Let $\alp=(\alpha^0,\ldots,\alpha^{r-1}) \in \mathcal{P}_{r,n}$ be
an $r$-partition of $n$ and let $\y^0,\ldots,\y^{r-1}$ be $r$ sets
of independent variables. We define
\begin{equation}\label{def-powsum}
p_{\alp}(\y^0,\ldots,\y^{r-1}):=\prodlim_{j=0}^{r-1}\prodlim_{i=1}^{\ell(\alpha^j)}\sumlim_{k=0}^{r-1}{\zeta^{j\cdot
k}p_{\alpha^{j}_{i}}(\y^k)},
\end{equation}
where $\zeta$ is a primitive $r^{\rm th}$ root of the unity. For example, for $r=4$ we get:
\begin{eqnarray*}
p_{(\alpha^0,\alpha^1,\alpha^2,\alpha^3)}(\y^0,\y^1,\y^2,\y^3) & =
&
\prodlim_{i=1}^{\ell(\alpha^0)}{(p_{\alpha^0_i}(\y^0)+p_{\alpha^0_i}(\y^1)+p_{\alpha^0_i}(\y^2)+p_{\alpha^0_i}(\y^3))}\\
& \cdot & \prodlim_{i=1}^{\ell(\alpha^1)}{(p_{\alpha^1_i}(\y^0)+\zeta p_{\alpha^1_i}(\y^1)+\zeta^2 p_{\alpha^1_i}(\y^2)+\zeta^3p_{\alpha^1_i}(\y^3))}\\
& \cdot & \prodlim_{i=1}^{\ell(\alpha^2)}{(p_{\alpha^2_i}(\y^0)+\zeta^2p_{\alpha^2_i}(\y^1)+\zeta^0p_{\alpha^2_i}(\y^2)+\zeta^2p_{\alpha^2_i}(\y^3))}\\
& \cdot &
\prodlim_{i=1}^{\ell(\alpha^3)}{(p_{\alpha^3_i}(\y^0)+\zeta^3p_{\alpha^3_i}(\y^1)+\zeta^2p_{\alpha^3_i}(\y^2)+\zeta
p_{\alpha^3_i}(\y^3))}.
\end{eqnarray*}

Now, for any $\lab$ and $\alp$ in $\mathcal{P}_{r,n}$, let $\chi
_{\alp}^{\lab}$ be the value of the irreducible character of $G$
indexed by $\lab$ on the conjugacy class of type $\alp$. The
following theorem presents the analogue of Frobenious formula for
$G$. Its proof can be deduced by the proof of formula 9.5''
in \cite[pp 177]{MD}. See also \cite{Sh}.

\begin{thm}\label{frob-G}

Let $\alp \in \mathcal{P}_{r,n}$ a $r$-partition of $n$. Then
$$p_{\alp}(\y^0,\ldots,\y^{r-1})=
\sumlim_{\lab \in \mathcal{P}_{r,n}}\chi _{\alp}^{\lab}s_{\la^0}
(\y^0)\cdots s_{\lambda^{r-1}}(\y^{r-1}).$$
\end{thm}
>From Theorem \ref{repsofgrpn} and Theorem \ref{frob-G} we obtain
the following identity.

\begin{thm}[Frobenius Formula for $H$] \label{frob-H}
Let $\alp \in \mathcal{P}_{r,n}$. Then
$$p_{\alp}(\y^0,\ldots,\y^{r-1})=\sumlim_{[\lab]} \left((\chi_{\alp}^{([\lab],\delta_1)}+\cdots +
\chi_{\alp}^{([\lab],\delta_{u(\lab)})}) \cdot \sumlim_{\vec{\mu}
\in [\lab]} {s_{\mu^0}(\y^0) \cdots
s_{\mu^{r-1}}(\y^{r-1})}\right),$$ where $[\lab]$ runs over all
the orbits of irreducible $G$-modules, $u(\lab)=|C_\lab|$, and $\muv$ runs over the
elements of the orbit $[\lab]$.
\end{thm}

\section{$n$-Orbital Standard Tableaux}

In this section we introduce a new class of standard Young
$r$-tableaux that will be useful for our purposes later on.

Let $\lab=(\la^0,\ldots,\la^{r-1})\in \mathcal{P}_{r,n}$ be an
$r$-partition of $n$. A {\em Ferrers diagram} of shape $\lab$ is
obtained by the union of the Ferrers diagrams of shapes
$\la^0,\ldots,\la^{r-1}$, where the $(i+1)^{\rm th}$ diagram lies
south west of the $i^{\rm th}$. A \emph{standard Young
$r$-tableau} $T:=(T^0,\ldots,T^{r-1})$ of shape $ \lab  $ is
obtained by inserting the integers $ 1,2,\ldots ,n $ as
\emph{entries} in the corresponding Ferrers diagram increasing
along rows and down columns of each diagram separately. We denote
by $ {\rm SYT}(\lab ) $ the set of all $r$-standard Young tableaux
of shape $ \lab $. Any entry in the $i$ component $T^i$ of $T \in {\rm SYT}(\lab )$ will be considered of color $i$.

A \emph{descent} in an $r$-standard Young tableau $ T $ is an
entry $ i $ such that $ i+1 $ is strictly below $ i $. We denote
the set of descents in $ T $ by $ \Des(T) $. Similarly to Section
\ref{flag} we let

\begin{eqnarray*}
d_{i}(T)&:=&|\{j\geq i\; :\; j\in \Des(T)\},\;  \; c_i=c_i(T):=k \; {\rm if} \; i \in T^k;\\
f_i(T)&:=& r \cdot d_i(T)+c_i(T), \;\; f(T):=(f_1(T),\ldots,f_n(T))\\
\col(T)&:=&c_1+\ldots +c_n,\\
\maj(T)&:=&\sum_{i\in\Des(T)}i, \; \;{\rm and} \;\; \fmaj(T):=r
\cdot \maj(T)+ \col(T).
\end{eqnarray*}

\noindent  For example, the tableau $T_1$ in Figure 1 belongs to
${\rm SYT}((1),(2),(2,1),(1,1),(3,1),(2))$. We have that  $
\Des(T)=\{1,3,5,8,11,12\} $, $ \maj(T)=40$, $\col(T)=1\cdot 2 +
2\cdot 3 + 3 \cdot 2 +4 \cdot 4 +5 \cdot 2=40$, and so
$\fmaj(T)=280$.

\begin{figure}\begin{center}\leavevmode
\epsfbox{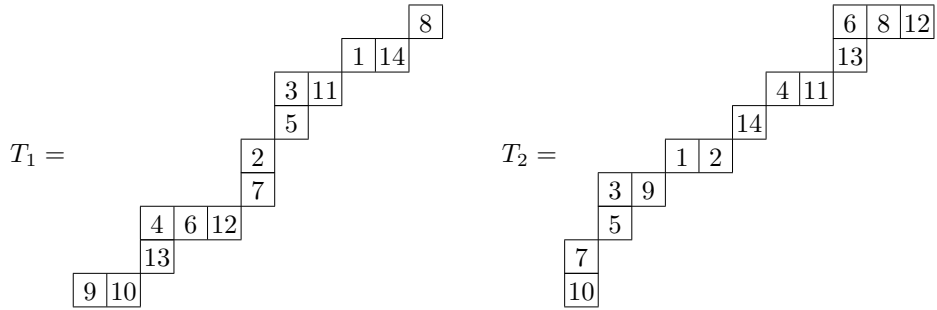}
\end{center}
\caption{Two $6$-standard Young tableaux}\label{1}
\end{figure}

Let $\lab=(\la^0,\ldots,\la^{r-1})\in \mathcal{P}_{r,n}$. As in Section \ref{repr-theory} let $[\lab]=\{\vec{\mu}\in \mathcal{P}_{r,n} \mid \vec{\mu} \sim \lab\}$ be the orbit of $\lab$ under the equivalence relation $\sim$ defined in (\ref{orbit}). An {\em orbital standard Young tableau} $T=(T^{0},\ldots,T^{r-1})$ {\em of type} $[\lab]$ is a standard Young $r$-tableau having one of the shapes
in $[\lab]$. The following definition is fundamental in our work. An {\em $n$-orbital standard Young tableau} of type $[\lab]$ is an orbital tableau of type  $[\lab]$ such that  $n \in T^0 \cup \cdots \cup T^{d-1}$. We denote by ${\rm OSYT}_n[\lab]$ the set of all $n$-orbital $r$-tableaux of type $\lab$.

More precisely, let $T$ be a standard Young $r$-tableau of shape $\lab$. From (\ref{orbit}),
it follows that all the possible orbital tableaux of type $[\lab]$, have shapes
obtained from that of $\lab$ by applying $i \cdot d$-times the shift operator (\ref{shift-op}) , for $i=0,\ldots,p-1$.
\begin{exa}{\rm Let $r=6$ and $n=17$. If $p=3$, and so $d=2$, then the two tableaux $T_1$ and $T_2$ in Figure 1 are of the same type $[\lab]=[(1),(2),(2,1),(1,1),(3,1),(2)]$: $T_1$ is $n$-orbital, while $T_2$ is not.
Differently, for $p=2$, and $d=3$ the two tableaux $T_1$ and $T_2$ are not in the same orbit $[\lab]$. Nervertheless, $T_2$ is an $n$-orbital tableau of type $[(3,1),(2),(1),(2),(2,1),(1,1)]$.}
\end{exa}


\section{Colored-descent representations}\label{col-des}

The module of coinvariants $\CC[\x]_H$ has a natural grading
induced from that of $\CC[\x]$. If we denote by $R_k$ its $k^{\rm
th}$ homogeneous component, we have:
\[\CC[\x]_H=\bigoplus_{k\geq 0} R_k.\]
Since the action (\ref{def-action}) preserves the degree, every homogeneous component $R_k$ is itself a $H$-module.
In this section we introduce a set of $H$-modules $R_{\mathcal D,
\mathcal C}$ which decompose $R_k$. The representations $R_{\mathcal D,
\mathcal C}$, called
{\em colored-descent representations}, generalize to all groups $G(r,p,n)$, the descent
representations introduced for $S_n$ and $B_n$ by Adin, Brenti and
Roichman in \cite{ABR}. See also \cite{BC2} for the case of $D_n$.
In the case of $G(r,n)$, a Solomon's descent algebra approach to these representations has been done by  Baumann and Hohlweg \cite{BH}. Since their study is
restricted to wreath products, it will be interesting to extend their
results to all complex reflection groups, thus getting
characters of all our modules as images of a particular class of elements of the group algebra.

We use most of the notions and tools introduced in Section
\ref{col-des-basis}. In particular, let see the first important feature of the colored-descent basis elements.

\begin{prop}
\label{stro} Let $\g \in \GA,$ and  $h \in H $. Then
\[ h \cdot \x_{\gamma}=\sum _{\{u\in
\GA \, :\, \lambda(\x_{u}) \unlhd \lambda (\x_{\gamma
})\}}n_{u}\x_{u}+y,\] where $ c_{u}\in \mathbb{C}$, and $ y \in
\mathcal{I}_{H} $.
\end{prop}
\begin{proof}
Expand $ M=h \cdot \x_{\g } $ as in (\ref{straightening}). Note
that $ f_{\mu (M^{\prime })}\not \in \mathcal{I}_{H} $ if and
only if $ \mu (M^{\prime })=\emptyset . $ Suppose that $ M^{\prime
} $ gives a nonzero contribution to this expansion of $ M $. This implies
 $ M^{\prime }=\x_{\gamma (M^{\prime })} $. If we
let $ u=\g (M^{\prime }) $, then $\la (\x_{u})=\la (M^{\prime
}) \unlhd \la (M)=\la (h \cdot \x_{\g })=\la (\x_{\g }).$
\end{proof}

If $|\la|=k$ then by Proposition \ref{stro}:
\begin{eqnarray*}
J_{\lambda }^{\unlhd} &:=&\textrm{span}_{\mathbb{C}}\{\x_{\gamma }+\mathcal{I}_{H}\, |\, \gamma
\in \GA,\, \lambda (\x_{\gamma })\unlhd \lambda \} \;\; {\rm and}\\
J_{\lambda }^{\lhd} &:=&\textrm{span}_{\mathbb{C}}\{\x_{\gamma }+\mathcal{I}_{H}\, |\, \gamma \in
\GA,\, \lambda (\x_{\gamma }) \lhd \lambda \}
\end{eqnarray*}
are submodules of $R_{k} $. Their quotient is still an
$H$-module, denoted by $R_{\lambda }:=\frac{J_{\lambda
}^{\unlhd}}{J_{\lambda }^{\lhd}}.$

For any $ \mathcal D\subseteq [n] $ we define the partition $\lambda_ {\mathcal D}:=(\lambda _{1},\ldots ,\la _{n}),$  where $ \la _{i}:=| \mathcal D\cap [i,n]|. $ For $\mathcal D \subseteq [n-1]$ and  $\mathcal C \in [0,r-1]^{n-1}\times[d]$, we define the vector
\[ \lambda _{\mathcal D, \mathcal
C}:=r \cdot \lambda _{ \mathcal D}+ \mathcal C,\] where sum stands
for sum of vectors. The following two observations are importantl in
our analysis.
\begin{rem}
\label{rem} For any $ \gamma, \gamma^{\prime} \in \GA $ we have:
\begin{itemize}
\item[$1)$] $\lambda(\x_{\gamma })=\lambda _{\mathcal D, \mathcal C},$  where $ \mathcal D=\Des(\gamma )$  and $\mathcal C=\Colrm(\gamma )$
\item[$2)$] $\lambda (\x_{\gamma })=\lambda (\x_{\gamma ^{\prime }})$
if and only if $\Des(\gamma )=\Des(\gamma ^{\prime })$ and $
\Colrm(\gamma )=\Colrm(\gamma ^{\prime }).$
\end{itemize}
\end{rem}



>From now on we denote
\[R_{\mathcal D, \mathcal C}:=R_{\la_{\mathcal D, \mathcal C}},\]
and by $ \bar{\x}_{\g}$  the image of the colored-descent basis element $ \x_{\g} \in
J^{\unlhd}_{\la_{\mathcal D, \mathcal C}}$ in the quotient $
R_{\mathcal D, \mathcal C} $. The next proposition follows from the definition of   $R_{\mathcal D, \mathcal C}$, and by Remark \ref{rem}.
\begin{prop}
\label{barra} For any $ \mathcal D \subseteq [n-1]$ and $
\mathcal C \in [0,r-1]^{n-1}\times[d] $, the set \[ \{\bar{\x}_{\g }\; :\; \g
\in \GA,\; \Des(\g )= \mathcal D\; {\rm and}\; \Colrm(\g)=
\mathcal C\}\] is a basis of $ R_{\mathcal D, \mathcal C} $.
\qed
\end{prop}

The $H$-modules $ R_{\mathcal D, \mathcal C} $ are called
\emph{colored-descent representations}. The dimension of $R_{\mathcal D, \mathcal C} $, is then given by the number of elements in $\Gamma$ with descent set  $\mathcal{D}$ and color vector $\mathcal{C}$. They decompose the $k^{\rm th}$ component
of $\CC[\x]_H$ as follows.

\begin{thm}
\label{decom} For every $ 0\leq k \leq r {{n}\choose{2}}+ n(d-1)
$,
\begin{equation}\label{isom}
R_{k}\cong \bigoplus _{\mathcal D, \mathcal C}R_{\mathcal D, \mathcal C}
\end{equation}
as $H$-modules, where the sum is over all $ \mathcal D\subseteq [n-1], \mathcal C \in [0,r-1]^{n-1}\times[d]$
such that \[ r\cdot \sum _{i\in \mathcal{D}}i+\sum_{j \in
\mathcal{C}} j =k.\]
\end{thm}
\begin{proof} By Theorem \ref{basis-coivariants}, the set $\{\x_\g +\mathcal{I}_H: \fmaj(\g)=k\}$ is a basis for $R_k$, and so the isomorphism (\ref{isom}) as $\CC$-vector spaces is clear. Now, let $ \lambda ^{(1)}<\lambda ^{(2)}<\cdots <\lambda
^{(t)} $ be a a linear extension of the dominance order on the
set $ \{\lambda (\x_{\gamma })\, :\, \fmaj(\gamma )=k\} $, i.e.,
$ \lambda ^{(i)}\lhd \lambda ^{(j)} $ implies $ i\leq
j $. Consider the following flag of subspaces of $ R_{k} $, $
M_{i}:=\textrm{span}_{\mathbb{C}}\{ \x_{\gamma }\, :\, \lambda
(\x_{\gamma })\leq \lambda ^{(i)}\}$ for $ i\in [r] $, and $
M_{0}:=0 $. By Proposition \ref{stro} $ M_{i} $ is a $ H
$-submodule of $ R_{k} $ and it is clear that $
M_{i}/M_{i-1}\cong R_{\lambda ^{(i)}} $ for $ i\in [r] $. By
Maschke's theorem we can easily prove, by induction on $ n-i $,
that\[ R_{k}\cong M_{i}\oplus R_{\lambda ^{(i+1)}}\oplus
R_{\lambda ^{(i+2)}}\oplus \cdots \oplus R_{\lambda ^{(t)}}\] for
all $ i\in [0,t] $. The claim follows for $ i=0 $, since $
M_{0}=0 $.
\end{proof}

\section{Stanley's generalized formula}\label{S-stanley}

In this section we prove a generalization of Stanley's formula \cite[Theorem 7.21.2]{StaEC2}, and a technical lemma that will be fundamental in the proof of Theorem \ref{main}. We need to introduce some more notation.

 A {\em reverse semi-standard Young tableau of shape} $\lambda$ is obtained by
 inserting positive integers into the diagram of $\lambda$ in such a
 way that the entries decrease weakly along rows and strictly down columns.
 The set of all such tableaux is denoted by ${\rm RSSYT(\la)}$.
 This notion generalizes to $r$-partitions as follows:
 Let $\lab=(\la^0,\ldots,\la^{r-1}) \in \mathcal{P}_{r,n}$.
 A {\em reverse semi-standard Young tableau of shape $\lab$} is an $r$-tuple of reverse semi-standard Young
 tableaux $\T=(\T^0,\ldots,\T^{r-1})$, where each $\T^i$ has shape $\la^i$, and every entry in $\T^i$ is congruent to $i+1$ (mod $r$), for all $0 \leq i \leq r-1$.
 The set of all such tableaux is denoted by ${\rm RSSYT}(\lab)$.
 To every $\T \in {\rm RSSYT}(\lab)$, we associate the {\em entries partition} of
 $\T$:
\begin{equation}\label{entry-p}
\theta(\T):=(\theta_1,\ldots,\theta_n),
\end{equation}
made of the entries of $\T$ in weakly decreasing order.

Similarly to \cite[Section 5]{ABR}, for every $\lab \in \mathcal{P}_{r,n}$ we define a map
\begin{eqnarray*}
\phi_\lab: {\rm RSSYT}(\lab) & \rightarrow &{\rm SYT}(\lab) \times \NN^n \\
\phi_\lab(\T)&\mapsto & (T,\Delta),
\end{eqnarray*}
as follows:
\begin{itemize}
\item[$T)$] Let $\theta(\T)=(\theta_1,\ldots,\theta_n)$ be the entries partition (\ref{entry-p}) of $\T$.
Then $T$ is the standard Young tableau of the same shape $\lab$ of
$\T$, with entry $i$ in the same cell in which $\T$ has
$\theta_i$, for all $i\in [n]$. If some of the entries of $\T$ are
equal then they are in different columns and the corresponding
entries of $T$ will be chosen increasing left to right.
\item[$\Delta)$] For every $i \in [n]$ let
\[\Delta_i:=\frac{\theta_i -f_i(T)- \theta_{i+1} + f_{i+1}(T)}{r},\]
and set $\theta_{n+1}=1$.
\end{itemize}

\begin{figure}\begin{center}\leavevmode
\epsfbox{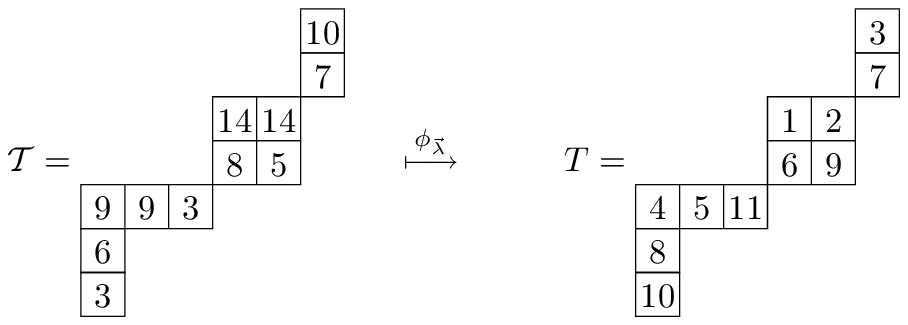}
\end{center}
\caption{The bijection $\phi_\lab$}\label{2}
\end{figure}

\begin{exa}{\rm  Let $r=3$, $n=11$ and let $\T=(\T^0,\T^1,\T^2)$ be the reverse standard
Young $3$-tableau of shape $\lab=((1,1),(2,2),(3,1,1))$ in Figure
\ref{2}. Note that each entry in $\T_i$ is congruent to $i+1$ (mod
3). The entries partition of $\T$ is
$\theta(\T)=(14,14,10,9,9,8,7,6,5,3,3)$. Let us compute the image of $\T$,
$\phi_\lab(\T)=(T,\Delta)$. The standard Young tableau $T$ is
drawn in Figure \ref{2}. Now, $\Des(T)=\{3,7,9\}$ and
$f(T)=(10,10,9,8,8,7,6,5,4,2,2)$. It follows that $(\theta(\T) -
f(T))_{i=1}^n=(4,4,1,1,1,1,1,1,1,1,1)$, and so
$\Delta=(0,1,0,0,0,0,0,0,0,0,0)$.}
\end{exa}

The following remark is not hard to see.
\begin{rem}\label{bijection-theta}
For every $\lab \in \mathcal{P}_{r,n}$, let $\phi_\lab(\T)=(T,\Delta)$. Then
\begin{itemize}
\item[$1)$] the map $\phi_\lab$ is a bijection;
\item[$2)$] $\theta_i -1=f_i(T)+ \sumlim_{j\geq i} r \Delta_j$ for all $i \in [n]$.
\end{itemize}
\end{rem}
\bigskip

Let $\y=y_1,y_2,\ldots$ be an infinite set of variables. Define $ {\mathbb {C}}[[y_{1},y_{1}y_{2},\ldots ]] $ to be the
ring of formal power series in the countably many variables $
y_{1},y_{1}y_{2},\ldots ,y_{1}y_{2}\cdots y_{k},\dots$.
A linear basis for $ {\mathbb {C}}[[y_{1},y_{1}y_{2},\ldots ]] $
consists of the monomials ${\bf y}^{\lambda }:=y_{1}^{\lambda
_{1}}\cdots y_{n}^{\lambda _{n}} $ for all partitions $ \lambda
=(\la _{1},\ldots ,\la _{n}) $.

Let $ \iota :{\mathbb {C}}[[y_{1},y_{1}y_{2},\ldots ]]\rightarrow
{\mathbb {C}}[[q_{1},q_{1}q_{2},\ldots ]] $ be the map defined on
the generators by \[ \iota ({\bf y}^{\la }):=\bf{q}^{\la ^{\prime }},\]
 and extended by linearity. Here as usual $\lambda^{\prime}$ denotes the conjugate partition of $\la$. Note that $ \iota  $ is not a ring
homomorphism.

The following Lemma is a generalization of Lemma 5.10 of
\cite{ABR}, which itself is a ``two dimensional" generalization of
Stanley's formula \cite[Theorem 7.21.2]{StaEC2}.

\begin{prop}[Stanley's generalized formula]\label{4.5}
Let $\lab \in \mathcal{P}_{r,n}$. Then
$$\iota [s_{\la^0}(\z^0) \cdot s_{\la^1}(\z^1) \cdots
s_{\la^{r-1}}(\z^{r-1})]= \sumlim_{T \in {\rm SYT}(\lab)}
\prodlim_{i=1}^{n}q_{i}^{f_{i}(T)}\frac{1}{{\prodlim_{i=1}^{n}(1-q_{1}^{r}q_{2}^{r}\cdots
q_{i}^{r})}},$$
where
\begin{eqnarray}
\z^0 & := & (1,y_1y_2 \ldots y_r,y_1 \ldots y_{2r},\ldots) \nonumber \\
\z^1 & := & (y_1,y_1y_2 \ldots y_{r+1},y_1y_2 \ldots y_{2r+1},\ldots) \label{name-z}\\
& \vdots & \nonumber \\
\z^{r-1} & := & (y_1y_2\ldots y_{r-1},y_1y_2\ldots
y_{2r-1},\ldots). \nonumber
\end{eqnarray}
\end{prop}
\begin{proof}
Let $\la$ be a partition of $n$. By \cite[Proposition
7.10.4]{StaEC2}, the Schur function $s_\la$ can be written as
\[s_\la(\y)=\sum_{\T \in {\rm RSSYT}(\la)} \prod_{i=1}^\infty y_i^{m_i(\T)},\]
where $m_i(\T)$ is the number of cells of $\T$ with entry $i$. Let
$\T=(\T^0,\ldots,T^{r-1})$ be a reverse semi-standard Young
$r$-tableaux of shape $\lab$. Note that for every
$i=0,\ldots,r-1$, the tableau $\frac{\T^i-(i+1)}{r}$ belongs to ${\rm
RSSYT}(\la^i)$, where $\frac{\T^i-(i+1)}{r}$ is obtained by
subtracting $i+1$ and then dividing by $r$ every entry of $\T^i$.

Now, let $\y^j=y_1^j,y_2^j,\ldots$, for $0 \leq j \leq r-1$, be a set of
$r$ sequences of independent variables. From the previous observation, it follows
\begin{equation}\label{mostro}
s_{\la^{0}}(\y^0)s_{\la^{1}}(\y^1) \cdots s_{\la^{r-1}}(\y^{r-1})=\sum_{\T}  \prod_{i=1}^\infty  (y^0_i)^{m_{ri-r+1}(\T^0)}{(y^1_i)}^{m_{ri-r+2}(\T^1)}\ldots{(y^{r-1}_i)}^{m_{ri}(\T^{r-1})},
\end{equation}
where $\T$ runs over all ${\rm RSSYT}(\lab)$. Now, if we
plug into (\ref{mostro}) the following substitutions:
\begin{center}
$\begin{array}{llll}
y^0_1= 1     &  y_2^0 = y_1\ldots y_{r}   & y_3^0 = y_1\ldots y_{2r} & \ldots \\
y^1_1= y_1 &  y_2^1= y_1\ldots y_{r+1}& y_3^1 = y_1\ldots y_{2r+1} & \ldots \\
 \vdots  & \vdots  &  \vdots  & \vdots \\
y^{r-1}_1  = y_1\cdots y_{r-1} &  y_2^{r-1}=y_1\ldots y_{2r-1} & y_3^{r-1} = y_1\ldots y_{3r-1} & \ldots
\end{array}$
\end{center}
we get
\[s_{\la^{0}}(\z^0)s_{\la^{1}}(\z^1) \cdots s_{\la^{r-1}}(\z^{r-1})=\sum_{\T \in {\rm RSSYT}(\lab)}  \prod_{i=1}^\infty y_i^{m_{> i}(\T)},\]
where $m_{> i}(\T)$ is the number of cells of $\T$ with entry strictly bigger then $i$, and $\z^i$ are defined as in (\ref{name-z}).
Let $\mu(\T)$ be the partition $(m_{> 1}(\T),m_{> 2}(\T),\ldots)$. Clearly $\mu_1(\T)\leq n$. It is no hard to see that the conjugate partition of $\mu(\T)$ is the entries partition of $\T$, with entries lowered by 1, i.e.,
\[\mu(\T)^{\prime}=(\theta_1-1,\ldots,\theta_n-1).\]
Now the result follows by Remark \ref{bijection-theta}.
\begin{eqnarray*}
\iota[s_{\la^{0}}(\z^0)s_{\la^{0}}(\z^1) \cdots
s_{\la^{r-1}}(\z^{r-1})]& = & \sum_{\T \in {\rm RSSYT}(\lab)}\prod_{i=1}^n q_i^{\theta_i-1}\\
& = & \sum_{\substack{T \in {\rm SYT}(\lab) \\(\Delta_1,\ldots,\Delta_n)\in \NN^{n}}}
\prod\limits_{i=1}^n q_i^{f_i(T)}\prod_{j=1}^n (q_1^r \cdots q_j^r)^{\Delta_j}.
\end{eqnarray*}
\end{proof}


\begin{lem}\label{4.6}
Let $T=(T^0,\ldots,T^{r-1})$ be a Young tableau of total size $n$
such that $n \not\in T^{r-1}$ and let $T^{\circlearrowleft}=(T^{r-1},T^0,\ldots,T^{r-2})$ be the Young
tableau obtained from $T$ by a $1$-shift. Then for all $i=1,\ldots,n$
$$f_i(T^{\circlearrowleft})=f_i(T)+1.$$
\end{lem}
\begin{proof}
We prove by backward induction on $n$. For $i=n$, since $n \not\in T^{r-1}$, we clearly have $f_n(T^{\circlearrowleft})=f_n(T)+1$. Hence let $i<n$.
If $i$ and $i+1$ are both in the same component, then we have
$c_i(T^{\circlearrowleft})-c_i(T)=c_{i+1}(T^{\circlearrowleft})-c_{i+1}(T).$
Now, if $i \in \Des(T)$ then also $i \in
\Des(T^{\circlearrowleft})$ thus $d_i(T^{\circlearrowleft})-d_i(T)=d_{i+1}(T^{\circlearrowleft})-d_i(T)$, and we conclude
$f_i(T^{\circlearrowleft})-f_i(T)=f_{i+1}(T^{\circlearrowleft})-f_{i+1}(T)=1$ by the induction hypothesis. If $i \not\in \Des(T)$ then the argument is similar.

We turn now to the more interesting case, where $i,i+1$ are not in
the same component. We assume first that $i,i+1 \not\in T^{r-1}$.
In this case, we clearly have
$c_i(T^{\circlearrowleft})-c_i(T)=c_{i+1}(T^{\circlearrowleft})-c_{i+1}(T).$ Since $i,i+1
\not\in T^{r-1}$, $i\in \Des(T)$ if and only if $i
\in \Des(T^{\circlearrowleft})$, and thus we have
$d_i(T^{\circlearrowleft})-d_i(T)=d_{i+1}(T^{\circlearrowleft})-d_{i+1}(T).$ Hence,
the result follows by the induction
hypothesis as above.

If $i \in T^{r-1}$ but $i+1 \not\in T^{r-1}$ then $i \not\in \Des(T)$ but $i
\in \Des(T^{\circlearrowleft})$ so we have
$d_i(T^{\circlearrowleft})-d_i(T)=d_{i+1}(T^{\circlearrowleft})+1-d_{i+1}(T).$
We also have $c_i(T^{\circlearrowleft})-c_i(T)=1-r$ and
$c_{i+1}(T^{\circlearrowleft})-c_{i+1}(T)=1$. Thus
$c_i(T^{\circlearrowleft})-c_i(T)=c_{i+1}(T^{\circlearrowleft})-c_{i+1}(T)-r.$
We finally have:
\begin{eqnarray*}
f_i(T^{\circlearrowleft})-f_i(T)& = & r\cdot(d_i(T^{\circlearrowleft})-d_i(T))+c_i(T^{\circlearrowleft})-c_i(T)\\
& = & r \cdot (d_{i+1}(T^{\circlearrowleft})+1-d_{i+1}(T))+c_{i+1}(T^{\circlearrowleft})-c_{i+1}(T)-r\\
& = & f_{i+1}(T^{\circlearrowleft})-f_{i+1}(T)=1
\end{eqnarray*}
by the induction hypothesis. The case $i \in T^r,i+1 \not\in T^r$ is similar.
\end{proof}

\section{Decomposition of $ R_{\mathcal D, \mathcal C} $}\label{S-decomposition}

In this section we prove a simple combinatorial description of the
multiplicities of the irreducible representations of $H$ in $ R_{\mathcal D, \mathcal C} $.

For an element $ h \in H$ let the (graded) trace of its action
on the polynomial ring $\CC[\x] $ be
\begin{equation}\label{trace-C}
{\rm Tr}_{\CC[\x]}(h):=\sum _{M}\langle h \cdot M,M\rangle \cdot {\bf q}^{\la (M)},
\end{equation}
where the sum is over all monomials $ M\in \CC[\x] $,
$ \la (M) $ is the exponent partition of $ M $, and the inner
product is such that the set of all monomials is an orthonormal basis
for $ \CC[\x] $. Note that $ \langle h \cdot M,M\rangle \in \{0, \zeta^{i} : i=0,\ldots,r-1\}$.

\begin{lem}\label{ciclo}
Let $h \in H$  be of cycle type
$\vec{\alpha}=(\alpha^0,\ldots,\alpha^{r-1})$. Then the graded trace of its action on
$\CC[\x]$ is
$${\rm Tr}_{\CC[\x]}(h)=\iota [p_{\vec{\alpha}}(\z^0,\z^1,\ldots,\z^{r-1})],$$ where the
parameters $\z^k$ are defined as in (\ref{name-z}).
\end{lem}
\begin{proof}
Let $h \in H$ be of cycle type $(\alpha^0,\ldots,\alpha^{r-1})$.
Decompose $h$ into cycles according to the cycle type of
$h$, as  explained in Section \ref{repr-theory}, to get
$h=\mathfrak{h}^0 \cdots \mathfrak{h}^{r-1}$, in such a way that $\mathfrak{h}^j$ is a product of cycles of type $j$, for each $0 \leq j \leq r-1$. We further decompose each $\mathfrak{h}^j$ into
disjoint cycles $\mathfrak{h}^{j,1},\ldots,\mathfrak{h}^{j,\ell(\alpha^j)}$.
For example, as we observed in Section \ref{repr-theory}, $h=62^5 4^4 3^1 1^6 5^3 \in G(8,2,6)$ can be decomposed as a product of $\mathfrak{h}^{1,1}=(1^6 6 5^3)$, $\mathfrak{h}^{5,1}=(3^1 4^4)$ and $\mathfrak{h}^{5,2}=(2^5)$.

Every monomial $M \in \CC[\x]$ satisfying $\langle h \cdot
 M,M \rangle \neq 0 $ is of the form $M=M^0 \cdots M^{r-1}$, where
the variables of $M^j$ have indices appearing in $\mathfrak{h}^j$ and all have same exponent. Once again, we can further decompose each $M^j$ as $M^j=M^{j,1} \cdots M^{j,\ell(\alpha^j)}$ according to the cycle decomposition of
$\mathfrak{h}^j$. It is easy to see that $h \cdot
M=\prodlim_{j=0}^{r-1}\prodlim_{i=1}^{\ell(\alpha^j)}{\mathfrak{h}^{j,i}\cdot
M^{j,i}}$. So, given such an $M$, we can calculate the contribution
of its parts to the trace separately.
Fix $0 \leq j \leq r-1$ and $1 \leq i \leq \ell(\alpha^j)$. To every $e^{j,i} \in \mathbb{N}$ corresponds
 a unique monomial $M^{j,i}$ with $\langle \mathfrak{h}^{j,i} \cdot M^{j,i} , M^{j,i} \rangle \neq
0$ such that $e^{j,i}$ is the common exponent for the factor $M^{j,i}$.
Moreover $\langle \mathfrak{h}^{j,i} \cdot M^{j,i},M^{j,i}\rangle=\zeta^{j\cdot e^{j,i}}$.
On the other hand, $\lambda(M^{j,i})$ consists of $\alpha_i^j$
copies of $e^{j,i}$, and we have
$$\iota^{-1}[{\bf q}^{\lambda(M^{j,i})}]={\bf y}^{\lambda(M^{j,i})'}=
{(y_1 \cdots y_{e^{j,i}})}^{\alpha^j_i}.$$ The last expression is
a summand in $p_{\alpha^j_i}(\z^{e^{j,i}})$. Here, $e^{j,i}$ is
taken mod $r$, and $\z^k$ is defined as in (\ref{name-z}).
Continuing our example, the monomial $M=x_1^2x_6^2x_5^2 \cdot
x_3^4x_4^4 \cdot x_2^5$ is such that $\langle h \cdot M,M \rangle
\neq 0 $. For $M^{5,1}=x_3^4x_4^4$, we have $\langle
\mathfrak{h}^{5,1} M^{5,1},M^{5,1} \rangle=\zeta^4$, and
$\iota^{-1}[{\bf q}^{\lambda(M_{5,1})}]={(y_1 \cdots y_4)}^2$,
with $(y_1\cdots y_4)^2$ appearing as a summand in $p_2(\z^4)$.

By summing over all the possibilities for $i,j$ and $e^{j,i}$,  we obtain
$${\rm Tr}_{\CC[\x]}(h)=
\iota\left( \prodlim_{j=0}^{r-1}\prodlim_{i=1}^{\ell(\alpha^j)}\sumlim_{k=0}^{r-1}\zeta^{j\cdot
k}p_{\alpha^j_i}(\z^k)\right)=\iota[p_{\vec\alpha}(\z^0,\z^1,\ldots,\z^{r-1})].$$
\end{proof}

Recall the notation introduced in Section \ref{col-des-basis}. Clearly, every monomial $N \in \CC[\x]$ can be written as
\begin{equation}\label{mono}
N=\vartheta_n^t\cdot M,
\end{equation}
with $M \in S:=\CC[\x]/(\vartheta_n)$, and $t \in \NN$.

>From the expansion (\ref{straightening}), it follows that the set
$\{\vartheta_n^t  \vartheta_{\mu}  \x_\gamma \}$ is an
homogeneous basis for $\CC[\x]$, where $\gamma \in \Gamma$, $\mu
\in \mathcal{P}(n-1)$ and $t \in \mathbb{N}$. Here
$\mathcal{P}(n-1)$ denotes the set of all partitions $\mu$ with
largest part at most $n-1$.

We now associate to any monomial $N=\vartheta_n^t \cdot M$ in
$\CC[\x]$ a triple
\begin{eqnarray}
\phi: \mathcal{M} & \longrightarrow & \Gamma(r,p,n)  \times \mathcal{P}(n-1)\times \NN \nonumber \\
N & \longmapsto & (\g(M), \mu(M), t).\label{bij}
\end{eqnarray}
Here $\mathcal{M}$ denote the set of all monomials in $\CC[\x]$, $\g(M)$
and $\mu(M)$ are the color index permutation and the complementary
partition of the monomial $M \in S$ defined in (\ref{mono}).
\begin{lem}
The map $\phi$ is a bijection
\end{lem}
\begin{proof}
The map is clearly into $\Gamma \times \mathcal{P}(n-1)  \times \NN$ since $\mu'(M)$ has at most $n-1$ parts.
On the other hand, to each triple $(\gamma,\mu,t) \in
\Gamma \times \mathcal{P}(n-1)  \times \NN$
we associate the monomial $N=\vartheta_n^t \cdot  M$, where $M$ is
the $\sqsubset$-maximal monomial in the expansion of $\vartheta_{\mu}
 \x_{\gamma}$. Clearly $\phi(N)=(\gamma,\mu,t)$.
\end{proof}

\begin{exa}
{\rm Consider the monomial  $N=x_1^{17}x_2^{14}x_3^7$. Then
$N=\vartheta_3^2\cdot M$, where $M=x_1^{11}x_2^{8}x_3$ is the
monomial appeared in Example \ref{ex-M}. It follows that
$\g(M)=1^52^23^1 \in \Gamma(6,2,3)$, and $\mu(M)=(2)$. Hence
$N\mapsto (1^52^2 3, (2),2)$.}
\end{exa}
Since $\lambda(N)=\lambda(M)+(dt)^n$ and
$\lambda(M)=\lambda(\x_{\gamma(M)}) + r\cdot \mu'(M)$, we get
\begin{equation}\label{exponent}
\lambda(N)=\lambda(\x_\g)+r \cdot \mu^\prime(M) + (dt)^n.
\end{equation}

Similarly to (\ref{trace-C}), we let
\begin{equation} \label{tra}
{\rm Tr}_{\CC[\x]_H}(h):=\sum _{\g \in \GA}\langle h \cdot (\x_{\g}+\mathcal{I}_{H}),\x_{\g }+\mathcal{I}_{H}\rangle \cdot {\bf q}^{\la (\x_{\g })},
\end{equation}
where the inner product is such that the colored-descent basis is orthonormal. The following result relates the two traces.
\begin{lem}
\label{ser} Let $ n \in \CC[\x]$ and $h \in H$. Then
\[{\rm Tr}_{\CC[\x]}(h)= {\rm Tr}_{\CC[\x]_H}(h) \cdot
\frac{1}{1-q_1^d \cdots q_n^d}
\prod_{i=1}^{n-1}{\frac{1}{1-q_1^r \cdots
q_i^r}}.\]
\end{lem}
\begin{proof}
We compute the trace (\ref{trace-C}) with respect to the
homogeneous basis $\{\vartheta_n^t  \vartheta_{\mu}
\x_\gamma \}$, (where $t \in \NN$, $\mu \in \mathcal{P}(n-1)$ and $\g \in \Gamma$),
using the inner product which makes this basis orthonormal.
 Note that by Proposition \ref{stro}, the action of $h$ on $\vartheta_n^t  \vartheta_{\mu}  \x_\gamma$ is triangular
 and thus the trace (\ref{trace-C}) now looks like
$${\rm Tr}_{\CC[\x]}(h)=\sum_{\substack{\mu \in \mathcal{P}(n-1) \\ \g \in \Gamma, \; t \in \NN} }{\langle h \cdot \vartheta_n^t
\vartheta_{\mu}  \x_{\g},\vartheta_n^t
\vartheta_{\mu}  \x_{\g} \rangle {\bf q}^{\lambda(N)}}$$ where $N$
is the $\sqsubset$-maximal monomial in the expansion of
$\vartheta_n^t
 \vartheta_{\mu}\x_\gamma$.

On the other hand, since $\vartheta_n^t \vartheta_{\mu}$ is $H$-invariant, we have $\langle h \cdot \vartheta_n^t \vartheta_{\mu} \x_{\g},\vartheta_n^t \vartheta_{\mu}
\x_{\g}\rangle=\langle h \cdot \x_{\g},\x_{\g}\rangle$, and thus from (\ref{tra}) and (\ref{exponent}) we get
\begin{eqnarray*}
{\rm Tr}_{\CC[\x]}(h)
& = &  {\rm Tr}_{\CC[\x]_H}(h)\cdot \sum_{\mu \in \mathcal{P}(n-1)}{\bf
q}^{r \cdot \mu'}\cdot
\sum_{t \in \NN}{{\bf q}^{(dt)^n}}\\
& = & {\rm Tr}_{\CC[\x]_H}(h) \cdot
\sum_{\ell(\la)\leq n-1}(q_1\cdots q_n)^{r\la}\cdot \sum_{t\geq
0}(q_1^d\cdots q_n^d)^t.
\end{eqnarray*}
\end{proof}




\begin{thm}\label{main}
For every $\mathcal{D} \subseteq [n-1]$ and $\mathcal{C} \subseteq [0,r-1]^{n-1}\times[d]$, $\lab
\in \mathcal{P}_{r,n}$ and $\delta \in C_{\lab}$, the
multiplicity of the irreducible representation of $\G$
corresponding to the pair $([\lab],\delta)$ in $R_{\mathcal{D},\mathcal{C}}$ is
$$|\{T \in {\rm OSYT}_n[\lab] \mid \Des(T)=\mathcal{D},\Colrm(T)=\mathcal{C}\}|.$$
\end{thm}
\begin{proof}
Let $h \in H$. If $h$ is of type $\alp=(\alpha^0,\ldots,\alpha^{r-1} )$ then
by Lemma \ref{ciclo} we have:

$${\rm Tr}_{\CC[\x]}(h)=\iota[p_{\alp}(\z^0,\ldots,\z^{r-1})].$$

By the Frobenius formula for $\G$, (Theorem \ref{frob-H}), we have

$$\iota[p_{\alp}(\z^0,\ldots,\z^{r-1})]=
\iota\left[\sumlim_{[\lab]}{(\chi_{\alp}^{([\lab],\delta_1)}+\cdots +
\chi_{\alp}^{([\lab],\delta_{u(\lab)})})} \cdot \sumlim_{\muv \in
[\lab]} {s_{\mu^0}(\z^0) \cdots s_{\mu^{r-1}}(\z^{r-1})}\right]$$
where $u(\lab)=|C_\lab|$,  $[\lab]$ runs over
the orbits of the irreducible $G$-modules, and
$\muv$ runs over the elements of the orbit of $\lab$.
Let us fix an irreducible $G$-module indexed by $\lab$ and
consider the expression $\sumlim_{\muv \in [\lab]}
{s_{\mu^0}(\z^0) \cdots s_{\mu^{r-1}}(\z^{r-1}})$. If we denote by
$\muv_0=(\mu_0^0,\ldots,\mu_0^r)$ one of the elements of $[\lab]$
and by $\delta_0^d$ the generator of $G/H$, then we can write
$[\lab]=\{(\muv_0)^{d i}) \mid 0 \leq i \leq p-1\}.$
By Lemma \ref{4.5} we have
\begin{equation}\label{mezza}
\iota[\sumlim_{\muv \in [\lab]} {s_{\mu^0}(\z^0) \cdots
s_{\mu^{r-1}}(\z^{r-1})}]=\sumlim_{\muv \in [\lab]}\sumlim_{T
\in {\rm SYT}(\muv)}\prodlim_{i=1}^{n}q_{i}^{f_{i}(T)}\frac{1}{{\prodlim_{i=1}^{n}(1-q_{1}^{r}q_{2}^{r}\cdots q_{i}^{r})}}.
\end{equation}

The action of $d$-shifts gives us a bijection between the sets
$\{T \in {\rm SYT}(\muv_0) \mid n \in T^0 \cup \cdots \cup
T^{d-1}\}$ and $\{T \in {\rm SYT}((\muv_0)^{\circlearrowleft d}) \mid n
\in T^{d} \cup \cdots \cup T^{2d-1}\}$. By iterating this
procedure we obtain that each set of the form $\{T \in {\rm
SYT}((\muv_0)^{\circlearrowleft di}) \mid n \in T^{id} \cup \cdots \cup
T^{(i+1)d-1}\}$ for $1 \leq i \leq p-1$ is in a bijective
correspondence with one of the sets  $\{T \in {\rm SYT}((\muv_0)^{\circlearrowleft dj})
\mid n \in T^0 \cup \cdots \cup T^{d-1}\}$ for some $0 \leq j \leq
p-1$.

Before finishing the proof, let us consider a
simple example. Let $G=G(6,n)$ and $H=G(6,3,n)$. Then the
generator of $G/H$ acts by a $2$-shift. The bijection described
above can be seen in the following diagram which has to be
considered as lying in a torus screen.
{\small $$
\xymatrix{\{T \in \Sno \mid n \in T^0 \cup T^1 \} \ar[dr]\ar[ddrr]
& \{T \in \Sno \mid n \in T^2 \cup T^3 \} & \{T \in \Sno \mid n
\in T^4 \cup T^5 \} &
\\
\{T \in \Stwo \mid n \in T^0 \cup T^1 \} \ar[dr] \ar[ddrr] & \{T
\in \Stwo \mid n \in T^2 \cup T^3 \} &
 \{T \in \Stwo \mid n \in T^4
\cup T^5 \}  &
\\
\{T \in \Sone \mid n \in T^0 \cup T^1 \} \ar[dr] \ar[ddrr] & \{T
\in \Sone \mid n \in T^2 \cup T^3 \} & \{T \in \Sone \mid n \in
T^4 \cup T^5 \} &
\\
\{T \in \Sno \mid n \in T^0 \cup T^1 \} \ar[dr] & \{T \in \Sno
\mid n \in T^2\cup T^3 \} & \{T \in \Sno \mid n \in T^4 \cup T^5
\}&
\\
\{T \in \Stwo \mid n \in T^0 \cup T^1 \} & \{T \in \Stwo \mid n
\in T^2 \cup T^3 \} & \{T \in \Stwo \mid n \in T^4 \cup T^5 \} }$$}
where we write for example $\Sno$ for ${\rm SYT}(\lambda^0,\ldots,\lambda^5)$.

By using Lemma \ref{4.6} together with the above argument we obtain that the RHS of (\ref{mezza}) is equal to
$$ \sumlim_{\muv \in [\lab]} \sumlim_{\substack{T \in
{\rm SYT}(\muv) \\ n \in T^0 \cup \cdots \cup T^{d-1}}}
\prodlim_{i=1}^{n}q_{i}^{f_{i}(T)}\frac{(1+{(q_1 \cdots q_n)}^d + {(q_1
\cdots q^n)}^{2d} + \cdots + {(q_1 \cdots q^n)}^{(p-1)d} )}{{\prodlim_{i=1}^{n}(1-q_{1}^{r}q_{2}^{r}\cdots
q_{i}^{r})}}.$$
It follows
\begin{eqnarray*}{\rm Tr}_{\CC[\x]}(\gamma)=
\frac{1-(q_1^d \cdots q_n^d)^p} {{(1-q_1^d \cdots q_n^d)
\prodlim_{i=1}^{n}(1-q_{1}^{r}q_{2}^{r}\cdots q_{i}^{r})}}
\sumlim_{[\lab]}{(\chi_{\alp}^{([\lab],\delta_1)}+\cdots +
\chi_{\alp}^{([\lab],\delta_{u(\lambda)})})} \\
\cdot \sumlim_{\muv
\in [\lab]}  \sumlim_{\substack{T \in {\rm SYT}(\muv) \\ n \in T^0 \cup
\cdots \cup T^{d-1}}} \prodlim_{i=1}^{n}q_{i}^{f_{i}(T)} .
\end{eqnarray*}
Finally, from Lemma \ref{ser} we obtain
\begin{eqnarray*}
{\rm Tr}_{\CC[\x]_H}(\gamma) &=&
\sumlim_{[\lab]}\left({(\chi_{\alp}^{([\lab],\delta_1)}+\cdots +
\chi_{\alp}^{([\lab],\delta_{u(\lambda)})})}
 \sumlim_{\muv
\in [\lab]}  \sumlim_{\substack{T \in {\rm SYT}(\muv) \\ n \in T^0 \cup
\cdots \cup T^{d-1}}} \prodlim_{i=1}^{n}q_{i}^{f_{i}(T)}\right)\\
&=& \sumlim_{[\lab]}{(\chi_{\alp}^{([\lab],\delta_1)}+\cdots +
\chi_{\alp}^{([\lab],\delta_{u(\lambda)})})} \sumlim_{T \in
{\rm OSYT}_n[\lab] } \prodlim_{i=1}^{n}q_{i}^{f_{i}(T)}.
\end{eqnarray*}

We conclude that the graded multiplicity of the irreducible $H$-module corresponding to the pair $([\lab],\delta)$ in $\CC[\x]_H$
is
$$\sumlim_{T \in{\rm OSYT}_n[\lab] } \prodlim_{i=1}^{n}q_{i}^{f_{i}(T)}=\sumlim_{T \in {\rm OSYT}_n[\lab]}
{\bf q}^{\lambda_{\Des(T),\Colrm(T)}}.$$
>From the proof of Theorem \ref{decom}, we obtain the decomposition
\[\CC[\x]_H\simeq \bigoplus_{\mathcal{D},\mathcal{C}}R_{\mathcal{D},\mathcal{C}},\]
as graded $H$-modules. By Proposition \ref{barra} and Remark \ref{rem} part $1)$, it follows that $R_{\mathcal{D},\mathcal{C}}$ is the homogeneous component of multidegree $\lambda_{\mathcal{D},\mathcal{C}}$ in $\CC[\x]_H$, and so we are done.
\end{proof}

As a consequence of Theorems \ref{decom} and \ref{main}, we obtain
the following result that was first proved by Stembridge \cite{Ste}), using a different terminology.
\begin{cor}
\label{stembrH} For $ 0\leq k\leq r {n \choose 2} + n(d-1)$,
the representation $ R_{k} $ is isomorphic to the direct sum
$ \oplus m_{k,(\lambda ,\delta)}V^{([\lab],\delta)} $,
where $ V^{([\lab],\delta)} $ is the irreducible representation
of $ H $ labeled by $([\lab],\delta)$, and \[
m_{k,([\lab],\delta)}:=\mid \{T\in {\rm OSYT}_n[\lab] \, :\, \fmaj(T)=k\}\mid .\]
\end{cor}

\section{Carlitz Identity}\label{identities}

In the case of classical Weyl groups and wreath
products, any major statistic is associated with a descent
statistic and their joint distribution is given by a nice closed
formula, called {\it Carlitz identity}. In this last section we
show that this is the case also for the complex reflection groups
$G(r,p,n)$.

 For any partition \( \la =(\la _{1},\ldots ,\la _{n})\in \mathcal{P}(n) \)
we let, for every \( j\geq 0 \),
\[m_{j}(\la ):=|\{i\in [n]\; :\; \la _{i}=j\}|,\;\; {\rm and} \;\;
{n\choose m_{0}(\la ),m_{1}(\la ),\ldots }\]
 be the multinomial coefficient.

\begin{lem}
Let $n \in \PP$. Then
\begin{equation}\label{pri}
\sum _{\ell (\la )\leq n}{n\choose m_{0}(\la ),m_{1}(\la ),\ldots }
\prod _{i=1}^{n}q_{i}^{\la _{i}}=\frac{\sumlim _{\g \in \Gamma}\prodlim _{i=1}^{n-1}q_{i}^{rd _{i}(\g )
+c_{i}(\g )}}{(1-q_{1}^d\cdots q_{n}^d)\prodlim _{i=1}^{n-1}(1-q_{1}^{r}\cdots q_{i}^{r})},
\end{equation}
in $\CC[[q_{1},\ldots ,q_{n}]]$.
\end{lem}

\begin{proof}
The LHS of the theorem is the multi graded Hilbert series of the
polynomials ring $\CC[{\bf x}]$ by exponent partition. In fact,
${n\choose m_{0}(\la ),m_{1}(\la ),\ldots }$ is the number of
monomials in $\CC[{\bf x}]$ with exponent partition equals to
$\la$. On the other hand, by the bijection (\ref{bij}), we have
\begin{eqnarray*}
\sumlim_{N \in \CC[\x]}{\bf q}^{\la(N)}& = &
\sum _{N=\vartheta_n^t M}{\bf q}^{\la (\x_{\g
(M)})+r\cdot \mu'(M)+(dt)^n}\\
& = & \sum _{\g \in \Gamma}{\bf q}^{\la
(\x_{\g })}\cdot \sum _{t\geq 0}(q_1^d\cdots q_n^d)^{t}\cdot \sum
_{\mu \in \mathcal{P}(n-1)}{q}^{r\mu^\prime} =
\frac{\sumlim _{\g \in \Gamma}\prodlim _{i=1}^{n-1}q_{i}^{rd
_{i}(\g ) +c_{i}(\g )}}{(1-q_{1}^d\cdots q_{n}^d)\prodlim
_{i=1}^{n-1}(1-q_{1}^{r}\cdots q_{i}^{r})}.
\end{eqnarray*}
\end{proof}

Recall the definition of flag-descent number given in (\ref{def-fdes}).
\begin{thm}[Carlitz identity for $G(r,p,n)$]\label{Ca-H} Let $n \in \NN$. Then
\begin{equation}\label{H-id}
\sum_{k \geq 0}[k+1]^n_q t^k=\frac{\sum_{h \in G(r,n,p)}
t^{\fdes(h)} q^{\fmaj(h)}}{(1-t)(1-t^rq^r)(1-t^rq^{2r})\cdots(1-t^rq^{(n-1)r})
(1-t^dq^{nd})}.
\end{equation}
\end{thm}
\begin{proof} The RHS is obtained by substituting $q_1=qt$, $q_2=\ldots=q_n=q$ in (\ref{pri}). The identity from the LHS of (\ref{H-id}) and the LHS of (\ref{pri}) is shown in \cite[Corollary 6.4]{ABR}.
\end{proof}
We refer to Theorem \ref{Ca-H} as the Carlitz identity for
$G(r,p,n)$. It is worth to note that the powers of the $q$'s in
the denominator, $r,2r,\ldots,(n-1)r,nd$, are actually the degrees
of $\G$.

Differently, if we plug $q_1=\ldots=q_n=q$ in (\ref{pri}), we get
the following identity

$$\frac{1}{(1-q)^n}=\frac{\sumlim_{\g\in
\Gamma}q^{\fmaj(\g)}}{(1-q^{nd})\prodlim_{i=1}^{n-1}(1-q^{ri})}.$$

Here, the LHS is the Hilbert series of the ring of polynomials
(simply graded) while the RHS is the product of the Hilbert series
of the module of coinvariants of $G(r,p,n)$ by the Hilbert series
of the invariants of $G(r,p,n)$ (on the denominator). Clearly, this equality reflects
the isomorphism between graded $H$-modules
$$\CC[\x] \simeq\CC[\x]^H \otimes  \CC[\x]_H,$$ which holds if and only if $H$ is a complex
reflection groups.

\appendix
\section{Appendix}

In this section we present another proof for the Carlitz identity
for $G(r,p,n)$, this time as a consequence of the Carlitz identity
for $G(r,n)$ which itself can be deduced either as a special case
of Theorem \ref{Ca-H} or by using a bijection analogues to
(\ref{bij}).  We start by presenting the Carlitz identity for
$G(r,n)$.

\begin{thm}[Carlitz identity for $G$]\label{Ca-G}
Let $n \in \NN$. Then
\begin{equation}
\sum_{k \geq 0}[k+1]^n_q t^k=\frac{\sum_{g \in \grn}
t^{\fdes(g)}q^{\fmaj(g)}}{(1-t)(1-t^rq^r)(1-t^rq^{2r})\cdots
(1-t^rq^{nr})}.
\end{equation}
\end{thm}

The proof we supply in this section points to an interesting
connection between the group $G(r,p,n)$ and its irreducible
representations. Actually, the idea of the proof of Theorem
\ref{main} about the decomposition of colored descent representations into irreducibles
arises here in the course of the alternative proof of Theorem
\ref{Ca-H}. More explicitly, the technical result we will use, Lemma
\ref{to-i}, carries the same role which Lemma \ref{4.6} carries for
standard Young tableaux. We start with some notations.

For every $0 \leq i \leq r-1$, we let
$$G_i=\{g=((c_1,\ldots,c_n),\sigma) \mid c_n=i\}.$$
Clearly we can decompose $G$ and $\GA$ as follows:
\begin{eqnarray}
G & = &G_0 \uplus G_1 \uplus \cdots \uplus G_{r-1}, \;\; {\rm and} \label{dec-G}\\
\GA& = &G_0 \uplus G_1 \uplus \cdots \uplus G_{d-1},\label{dec-H}
\end{eqnarray}
where $\uplus$ stands for disjoint union.

For $g=((c_1,\ldots,c_n),\s) \in G$ and $i \in \NN$ define
\[g^i:=((\tilde{c}_1,\ldots,\tilde{c}_n),\s),\]
where $\tilde{c}_k \equiv c_k+i$ mod $r$.

Moreover, for any $g=((c_1,\ldots,c_n),\s) \in \grn$ we let $\h=((c_1,\ldots,c_{n-1}),\hat{\s}) \in G(r,n-1)$, where for all $i \in [n-1]$
\[ \hat{\s}(i):=\left\{ \begin{array}{ll} \s(i) & {\rm if} \;\; \s(i)<\s(n)\\
                                                     \s(i)-1 & {\rm if} \;\; \s(i)>\s(n). \end{array}\right.
                                                     \]
For example, let $g=((4,1,3,0,2,1),416253) \in G(5,6)$. Then $g^2=((1,3,0,2,4,3),416253) \in G(5,6)$, and $\h=((4,1,3,0,2),31524) \in G(5,5)$. It is easy to see that
\begin{equation}\label{feature}
\Des(\h)=\Des(g) \cap [n-2] \;\; {\rm and} \;\; \Colrm(\h)=\Colrm(g) \cap [n-1].
\end{equation}

\begin{lem}\label{to-i}
Let $g \in G_{0}$. Then for every $i<r$:
\begin{eqnarray}
\fdes(g^i) & = & \fdes(g)+i,\\
\fmaj(g^i) & = & \fmaj(g)+ni.
\end{eqnarray}

\end{lem}

\begin{proof}
We proceed by induction on $n$. If $n=1$ then the statement is
clearly true. Hence suppose $n>1$, and let
$g=((c_1,\ldots,c_n),\s)$. There are three cases to consider.
\begin{enumerate}
  \item [1)] $\s(n-1)>\s(n)$.
\end{enumerate}
In this case, $c_{n-1}=c_n=0$, and so  $\h \in G_0(r,n-1)$. From
(\ref{feature}) we have:
\[ \Des(g)=\Des(\h)\cup \{n-1\} \;\; {\rm and} \;\;  \Des(g^i)=\Des(\h^i)\cup
\{n-1\}.\] Hence  $\des(g)=\des(\h)+1$, $\maj(g)=\maj(\h)+n-1$,
and analogously $\des(g^i)=\des(\h^i)+1$,
$\maj(g^i)=\maj(\h^i)+n-1$. Since $c_1(g)=c_1(\hat{g})$ and
$c_1(g^i)=c_1(\hat{g}^i)$, by induction we get:
\begin{eqnarray*}
\fdes(g^i) & = & r \cdot \des(g^i)+c_1(g^i) = r \cdot \des(\h^i)+r +c_1(\h^i)\\
            & = & \fdes(\h^i)+r = \fdes(\h)+i+r \\
            & = & \fdes(g)+i.
\end{eqnarray*}
On the other hand, from (\ref{feature}) we get:
$\col(g^i)=\col(\h^i)$. It follows
 that: \begin{eqnarray*}
\fmaj(g^i) & = & r \cdot \maj(\h^i)+rn-r +\col(\h^i)\\
            & = & \fmaj(\h^i)+rn-r = \fmaj(\h)+in+rn-r \\
            & = & \fmaj(g)+in.
\end{eqnarray*}
\begin{enumerate}
  \item [2)] $\s(n-1)<\s(n)$ and $c_{n-1}=c_n$.
\end{enumerate}
In this case $\Des(g)=\Des(\h)$ and  $\Des(g^i)=\Des(\h^i)$. Since
$\h \in G_0(r,n-1)$, the result easily follows by induction just
as in case 1).
\begin{enumerate}
  \item [3)] $\s(n-1)<\s(n)$ and $c_{n-1}=k\neq 0$.
\end{enumerate}
Suppose first that $i=r-k$. In this case the position $n-1$, which
is not a descent for $g$, becomes a descent for $g^i$, hence
$\Des(g^i)=\Des(\h^i)\cup \{n-1\}.$ It follows that:
\begin{equation}\label{becomedes}
\fdes(g)=\fdes(\h) \;\; {\rm and} \;\; \fdes(g^i)=\fdes(\h^i)+r.
\end{equation}
Now, $c_{n-1}(\h^i)=c_{n-1}+r-k=0$, hence we can apply induction.
Then for all $j < r$ we have that
$\fdes(\h^{r-k+j})=\fdes(\h^{r-k})+j$. In particular if we set
$j=k$ then we obtain:
\begin{equation}\label{giro}
\fdes(\h)=\fdes(\h^{r-k+k})=\fdes(\h^{r-k})+k.
\end{equation}
Hence, from (\ref{becomedes}) and (\ref{giro}) we obtain:
\[\fdes(g^{r-k})=\fdes(\h^{r-k})+r=\fdes(\h)+r-k=\fdes(g)+r-k.\]
On the other hand it is clear that
$\maj(g^{r-k})=\maj(\h^{r-k})+n-1$ and that
$\col(g^{r-k})=\col(\h^{r-k})+r-k$. Hence
\begin{eqnarray}
\fmaj(g^{r-k}) & = & r \cdot \maj(g^{r-k})+\col(g^{r-k}) \nonumber \\
                & = & r \cdot \maj(\h^{r-k})+ rn -r +
                \col(\h^{r-k})+r-k \nonumber \\
                & = & \fmaj(\h^{r-k})+rn-k.  \label{sopra}
\end{eqnarray}
Now, since $\h^{r-k} \in G_0(r,n-1)$, similarly to (\ref{giro})
we obtain \begin{equation}\label{uguale}
\fmaj(\h)=\fmaj(\h^{r-k})+k(n-1).
\end{equation}
Since $\fmaj(\h)=\fmaj(g)$, the result follows from (\ref{sopra})
and (\ref{uguale}).The other two cases: $i
\lessgtr r-k$ are similar and are thus left to the reader.
\end{proof}

\begin{proof}[{\bf Second proof of Theorem \ref{Ca-H}}]
By Lemma \ref{to-i}, and the decompositions (\ref{dec-G}) and (\ref{dec-H}) we obtain:
\begin{eqnarray*}
\sum_{g \in \grn}t^{\fdes(g)}q^{\fmaj(g)}& = &\sum_{g \in
G_0}t^{\fdes(g)}q^{\fmaj(g)}(1 + tq^{n}+ t^2q^{2n}+\ldots
+t^{r-1}q^{(r-1)n}), \;\; {\rm and}\\
\sum_{g \in G(r,p,n)}t^{\fdes(g)}q^{\fmaj(g)}& =& \sum_{g \in G_0}t^{\fdes(g)}q^{\fmaj(g)} (1 + tq^{n}+t^{2}q^{2n}+\ldots +t^{(d-1)}q^{(d-1)n}).
\end{eqnarray*}
It follows now that:
\[\sum_{g \in \grn}t^{\fdes(g)}q^{\fmaj(g)}=
\sum_{g \in G(r,p,n)}t^{\fdes(g)}q^{\fmaj(g)} (1 + t^{d}q^{dn}+
t^{2d}q^{2dn}+\ldots +t^{(p-1)d}q^{(p-1)dn}),\]
and Theorem \ref{Ca-G} ends the proof.
\end{proof}

\noindent {\bf Acknowledgments.} The first author wants to thank Ron Adin, Yuval Roichman, and Mishael Sklarz for useful discussions, and
Alex Lubotzky and the Einstein Institute of Mathematics
at the Hebrew University for hosting his stay.
The second author would like to thank Christophe Hohlweg for some helpful discussions about Solomon descent algebra.


\bigskip

\begin{tabular}{lll}
Eli Bagno &  \quad \quad \quad & Riccardo Biagioli \\
Einstein Institute of Mathematics &  & LaCIM \\
The Hebrew University &  & Universit\'e du Qu\'ebec \`a Mont\-r\'eal\\
Givaat Ram, Jerusalem, ISRAEL &  & Montr\'eal, Qu\'ebec, H3C 3P8, CANADA\\
email: $\mathtt{bagnoe@math.huji.ac.il}$ &  & email: $\mathtt{biagioli@math.uqam.ca}$
\end{tabular}
\bigskip

\end{document}